\documentclass[12pt]{amsart}

\usepackage{amssymb}

\usepackage{enumerate}

\makeatletter
\@namedef{subjclassname@2010}{%
  \textup{2010} Mathematics Subject Classification}
\makeatother

\newtheorem{thm}{Theorem}[section]
\newtheorem{cor}[thm]{Corollary}
\newtheorem{lem}[thm]{Lemma}

\theoremstyle{definition}

\numberwithin{equation}{section}

\begin{document}





\title[Baker type lower bounds]{On Baker type lower bounds for linear forms}

\author[T. Matala-aho]{Tapani Matala-aho}
\address{Matematiikan laitos\\
Pl 3000, 90014 Oulun Yliopisto\\
 Finland}
\email{tapani.matala-aho@oulu.fi}

\date{}

\begin{abstract}
A criterion is given for studying (explicit) Baker type lower bounds of linear forms
in numbers $1,\Theta_1,...,\Theta_m\in\mathbb{C}^*$
over the ring $\mathbb{Z}_{\mathbb{I}}$ of an imaginary quadratic field $\mathbb{I}$.
This work deals with the simultaneous auxiliary functions case.
\end{abstract}



\maketitle

\section{Introduction}

We give a criterion for studying (explicit) Baker type lower bounds of linear forms
in given numbers $\Theta_0,...,\Theta_m\in\mathbb{C}^*$.
Throughout this work, let $\mathbb{I}$ denote an imaginary quadratic field with
$\mathbb{Z}_{\mathbb{I}}$  it's ring of integers.
By an explicit Baker type lower bound we mean any positive lower bound
\begin{equation}\label{Bakertype1}
\left|\beta_0 \Theta_0+...+\beta_m \Theta_m\right| > F(H_0,...,H_m,m)
\end{equation}
valid for all 
$\overline{\beta}=(\beta_0,...,\beta_m)^T\in\mathbb{Z}_{\mathbb{I}}^{m+1}\setminus\{\overline{0}\}$
with
$\prod_{j=0}^{m} H_j\ge \hat H\ge 1,\quad H_j\ge h_j= \max\{1,|\beta_j|\},$
where the dependence on each individual term $H_0,...,H_m$, $m$ and numbers $\Theta_0,...,\Theta_m$ is
explicitly given in the functional dependence $F(H_0,...,H_m,m)$ and 
the dependence on $\Theta_0,...,\Theta_m,m$ is explicitly given in the constant 
$\hat H=\hat H(\Theta_0,...,\Theta_m,m)$.

With the assumption that $\gamma_0, \gamma_1,...,\gamma_m\in\mathbb{Q}^*$ are distinct,
Baker  \cite{BAKER} proved that there exist  positive constants $\delta_1,\delta_2$ and $\delta_3$ such that
\begin{equation}\label{Bakerresult}
\left|\beta_0 e^{\gamma_0}+...+\beta_m e^{\gamma_m}\right| > 
\frac{\delta_1 M^{1-\delta(M)}}{\prod_{j=0}^{m} h_j},
\end{equation}
for all 
$\overline{\beta}=(\beta_0,...,\beta_m)^T \in \mathbb{Z}^{m}\setminus\{\overline{0}\},
\quad h_j= \max\{1,|\beta_j|\},$
with
\begin{equation}\label{Bepsilon}
\delta(M)\le \frac{\delta_2}{\sqrt{\log\log M}},\quad 
M=\underset{0\le j\le m}\max\{|\beta_j|\}\ge \delta_3>e.
\end{equation}
Here we note that  the constants $\delta_1,\delta_2,\delta_3$ in Baker's work \cite{BAKER}
are not explicitly given. Mahler  \cite{MAHLER} made Baker's result completely explicit.

There are many subsequent works, where the authors prove Baker type lower bounds for
values of functions belonging usually to a class of
Siegel's $E$- or $G$-functions or $q$-hypergeometric functions evaluated at rational points, 
see e.g. \cite{FELDMAN1}, \cite{FELDMAN3}, \cite{VAAZUD} and \cite{ZUDILIN}.
For a more comprehensive list of references, see \cite{FELDMAN3}.
In the above mentioned works Siegel's lemma is a standard tool for producing 
a first or second kind Pad\'e-approximation construction of certain auxiliary functions.
These constructions correspond to one linear form (one auxiliary function)
or simultaneous linear forms (several auxiliary functions). 

In this work we shall not do such constructions but we are interested in the next step.
Namely, how to use appropriate linear forms to prove Baker type lower bounds?
We shall answer the above question by giving a criterion
in the simultaneous linear forms case.

Let us descripe our criterion in a nutshell.
Fix $\Theta_1,...,\Theta_m\in\mathbb{C}^*$ and put
$\ \overline{n}=(n_1,...,n_m)^T,\ N=N(\overline{n})=n_1+...+n_m.$
Assume that we have a sequence of simultaneous linear forms
\begin{equation}\label{SIMLINFORMS}
L_{k,j}(\overline{n})=A_{k,0}(\overline{n})\Theta_j+A_{k,j}(\overline{n}),
\quad k=0,1,...m,\ j=1,...,m,\quad \overline{n}\in\mathbb{Z}_{\ge 1}^{m},
\end{equation}
where
$A_{k,j}=A_{k,j}(\overline{n})\in\mathbb{Z}_{\mathbb{I}}$
satisfy a certain determinant condition. Suppose also that
\begin{equation}\label{ESTIMATE1}
|A_{k,0}(\overline{n})|\le e^{(aN+b\log N)g(N)+b_0N(\log N)^{1/2}+b_1N+b_2\log N+b_3},
\end{equation}
\begin{equation}\label{ESTIMATE2}
|L_{k,j}(\overline{n})|\le e^{(dN-cn_j)g(N)+ e_0N(\log N)^{1/2}+e_1N+e_2\log N+e_3},
\end{equation}
for $k,j=0,1,...,m$, where $a,b,c,d,b_i,e_i$ are non-negative parameters satisfying
$a, c-dm>0$.
Then, in the cases $g(N)\in\{1,\log N, N\}$, we shall 
prove that there exist explicit positive constants $F_l, G_l$ ($l\in\{1,2,3\}$), such that 
\begin{equation}\label{BTLB1}
\left|\beta_0+\beta_1\Theta_1+...+\beta_m\Theta_m\right| > 
F_l \left(\prod_{j=1}^{m}(2mH_j)\right)^{-\frac{a}{c-dm}-\epsilon_l(H)}
\end{equation}
holds for all 
$\overline{\beta}=(\beta_0,\beta_1,...,\beta_m)^T \in 
\mathbb{Z}_{\mathbb{I}}^{m+1}\setminus\{\overline{0}\}$
and
$
H=\prod_{j=1}^{m}(2mH_j)\ge G_l,\quad  H_j\ge h_j=\max\{1,|\beta_j|\} 
$
with an error term
$\epsilon_l(H)\underset{H\to\infty}\to 0.$
The constants $F_l, G_l$ and the error term will be given explicitly in terms of the parameters 
$a,b,c,d,b_i,e_i$ and, in particular of $m$.

The underlying idea behind our treatment is well known already from Baker's work \cite{BAKER}.
Namely, the idea, see formula (22) in \cite{BAKER}, is to fix the parameter $n_j$ with 
the corresponding  individual height $H_j$ (in our notation). 
In our work we shall express this phenomenon first in a nutshell, see \eqref{ONEPAGE},
and then in a refined form, see \eqref{Ssim}.

An advantage of our treatment compared with existing treatments
is that one may easily see if the contribution to the lower
bound is coming from the Diophantine method itself or from the auxiliary construction.
For example, apart from the condition $n_1+...+n_m=N$, we don't need 
any extra condition between $n_j$ and $N$ in our treatment.
Of course, some extra conditions may be needed for good auxiliary constructions.
In particular, this is the case when Siegel's lemma is involved.
See, e.g. \cite{VAAZUD}, formula (14), where the authors additionally assume that 
$n_j>\delta N$, $j=1,...,m$, for some  $0<\delta<1/m$.
In \cite{SANKILAMPI} the corresponding condition reads $n_j>2N/(\log N)$, $j=1,...,m$, 
formula (4) in Chapter III.
In \cite{AKT}, however, you may find a slightly different approach.

Our Theorems \ref{BAKERAXIOM1}, \ref{BAKERAXIOM2} and \ref{BAKERAXIOM3} are designed to be applied 
in the following manner. 
Let $f(z)$ be a $G$-, $E$- or $q$-hypergeometric function and denote $\Theta_1=f(\alpha_1),...,\Theta_m=f(\alpha_m)$,
$\alpha_1,...,\alpha_m\in\mathbb{I}^*$.
Suppose that one can construct simultaneous linear forms of the type \eqref{SIMLINFORMS} satisfying the estimates
\eqref{ESTIMATE1} and \eqref{ESTIMATE2} with a certain determinant condition, then 
our Theorem \ref{BAKERAXIOM1}, \ref{BAKERAXIOM2} or \ref{BAKERAXIOM3} will give a corresponding Baker type
lower bound \eqref{BTLB1}.
So far our results (Theorems \ref{BAKERAXIOM2} and \ref{BAKERAXIOM3})
have been applied in the works \cite{AKT} and \cite{LEINONEN}.
In \cite{AKT}, Ernvall-Hyt\"onen, Lepp\"al\"a and Matala-aho 
constructed simultaneous linear forms of the type \eqref{SIMLINFORMS}
(satisfying conditions \eqref{ESTIMATE1}-\eqref{ESTIMATE2} with $g(N)=\log N$)
for the exponential function values $e^{\alpha_0}, e^{\alpha_1}, ..., e^{\alpha_m}$,
where $\alpha_0,..., \alpha_m\in \mathbb{I}$. 
(Note that the exponential function belongs to the class of Siegel's $E$-functions.) 
By applying Theorem \ref{BAKERAXIOM2} of the present paper the authors in \cite{AKT} 
proved substantial improvements
of  the explicit versions, see Mahler \cite{MAHLER} and Sankilampi \cite{SANKILAMPI},
of Baker's work \cite{BAKER} about exponential values at rational points. 
In particular, the dependence on $m$ is improved.
As an example from the work \cite{AKT} we mention a new explicit Baker type lower bound
\begin{equation*}
\left|\beta_0+\beta_1 e+\beta_1 e^{2}+...+\beta_me^{m}\right| >
\frac{1}{h^{1+\hat \epsilon(h)}},\quad  h=h_1\cdots h_m,
\end{equation*}
valid for all
$\overline{\beta}=(\beta_0,...,\beta_m)^T \in \mathbb{Z}_{\mathbb{I}}^{m}\setminus\{\overline{0}\},
\ h_i= \max\{1,|\beta_i|\}$ 
with
\begin{equation*}
\hat \epsilon(h)=\frac{(4+7m)\sqrt{\log(m+1)}}{\sqrt{\log\log h}},\ 
\log h \ge m^2(41\log(m+1)+10)e^{m^2(81\log(m+1)+20)}.
\end{equation*}
As far as we know, the published dependences on $m$ in  $\hat \epsilon(h)$ have been at least quadratic
and in lower bound of $\log\log h$ at least quartic.
The second application of our work is presented in Leinonen's work \cite{LEINONEN}.
In a pioneer work \cite{VAAZUD} V\"a\"an\"anen and Zudilin proved Baker type
results for a class of $q$-hypergeometric series. 
Following the work \cite{VAAZUD}, Leinonen \cite{LEINONEN} constructed simultaneous 
linear forms of the type \eqref{SIMLINFORMS}
(satisfying conditions \eqref{ESTIMATE1}-\eqref{ESTIMATE2} with $g(N)=N$) 
and proved some generalizations of the results in \cite{VAAZUD}.
Moreover, in \cite{LEINONEN} Leinonen applied our Theorem \ref{BAKERAXIOM3} 
with her linear forms and gave explicit Baker type lower bounds which
sharpened her results as well the results of V\"a\"an\"anen and Zudilin.

\section{Background from metrical theory}

From the general metrical theory, see \cite{BAKER2},
\cite{CASSELS}, \cite{FELDMAN3}, \cite{SCHMIDT}, \cite{SHIDLOVSKII}
we get the following well known results.

\begin{thm}
Let $1,\Theta_1,...,\Theta_m\in\mathbb{R}$ be linearly independent over $\mathbb{Q}$. 
Then there exist infinitely many primitive vectors
$(\beta_0,...,\beta_m)^T\in\mathbb{Z}^{m+1}\setminus\{\overline{0}\}$
with 
$h_j:=\max\{1,|\beta_j|\},\ j=1,...,m,$
satisfying
\begin{equation*}\label{}
\left|\beta_0+\beta_1 \Theta_1+...+\beta_m \Theta_m\right|<\frac{1}{\prod_{j=1}^{m} h_j}.
\end{equation*}
\end{thm}

In the complex case
Shidlovskii \cite{SHIDLOVSKII} studies linear forms over the ring of 
rational integers and gives the following result.
\begin{thm}\cite{SHIDLOVSKII}\label{SHIDLOVSKII}
Let
$\Theta_0=1,\Theta_1,...,\Theta_m \in\mathbb{C}$
and
$H\in \mathbb{Z}_{\ge 1}$
be given. Then there exists a non-zero rational integer vector
$(\beta_0,\beta_1,...,\beta_m)^T\in\mathbb{Z}^{m+1}\setminus\{\overline{0}\}$
with
$|\beta_j|\le H$, $j=0,1,...,m,$
satisfying
\begin{equation*}\label{}
|\beta_0+\beta_1\Theta_1+...+\beta_m\Theta_m|
\le \frac{c}{H^{(m-1)/2}},\quad c=\sqrt{2}\sum_{j=0}^{m}|\Theta_j|.
\end{equation*}
\end{thm}

We are interested in linear forms over the ring of integers $\mathbb{Z}_{\mathbb{I}}$
in an imaginary quadratic field $\mathbb{Q}(\sqrt{-D})$, $D\in\mathbb{Z}_{\ge 1}$, 
$D\not\equiv 0\pmod 4$.
For that purpose we prove
\begin{thm}\label{ZC}
Let 
$\Theta_1,...,\Theta_m \in\mathbb{C}$
and
$H_1,...,H_m \in \mathbb{Z}_{\ge 1}$
be given. Then there exists a non-zero integer vector
$
(\beta_0,\beta_1,...,\beta_m)^T\in\mathbb{Z}_{\mathbb{I}}^{m+1}\setminus\{\overline{0}\}
$
with
$|\beta_j|\le H_j$, $j=1,...,m$,
satisfying
\begin{equation}\label{linform1C}
|\beta_0+\beta_1\Theta_1+...+\beta_m\Theta_m|
\le \left(\frac{2^{\tau}D^{1/4}}{\sqrt{\pi}}\right)^{m+1}\frac{1}{H_1\cdots H_m},
\end{equation}
where
$\tau=1$, if $D\equiv 1$ or $2\pmod 4$ and $\tau=1/2$, if $D\equiv 3\pmod 4$.
\end{thm}

\section{Results}

\subsection{A general target}

Let $f(z)$ belong to one of the following classes of functions, enumerated by 1--3.\\
1. The class of Siegel's $G$-functions. 
Typical examples are logarithm and Gauss hypergeometric functions and more generally
non-entire hypergeometric series.\\
2. The class of Siegel's $E$-functions. 
Typical examples are exponential and Bessel functions and more generally entire
hypergeometric series.\\
For definition of Siegel's $E$- and $G$-functions we refer to \cite{FELDMAN3}.\\
3.  The $q$-hypergeometric series.
Typical examples are\\ 
$\sum\limits_{n=0}^{\infty}q^{n^2}$ and  
$\sum\limits_{n=1}^{\infty}1/\prod\limits_{i=1}^{n}(1-q^i)$, $|q|<1$.\\

Our Theorems \ref{BAKERAXIOM1}, \ref{BAKERAXIOM2} and \ref{BAKERAXIOM3} are designed to be applied 
in the following manner. 
Denote $\Theta_1=f(\alpha_1),...,\Theta_m=f(\alpha_m)$, $\alpha_1,...,\alpha_m\in\mathbb{I}^*$.
Suppose that one can construct simultaneous linear forms of the type \eqref{LINFORMS} satisfying the conditions
\eqref{DET}, \eqref{accdm}, \eqref{AeqN} and \eqref{Lerin}, then 
our Theorem \ref{BAKERAXIOM1}, \ref{BAKERAXIOM2} or \ref{BAKERAXIOM3} will give  
a Baker type lower bound for the quantity
\begin{equation}\label{Bakertype12}
\left|\beta_0+\beta_1 \Theta_1+...+\beta_m \Theta_m\right|. 
\end{equation}
It is a general phenomenon in the field of Diophantine approximations that 
Pad\'e approximations and Siegel's lemma give estimates of shape \eqref{AeqN} and \eqref{Lerin}.
However, often it is hard to find such bounds if the condition \eqref{accdm} holds, too.

\subsection{A criterion}

Fix now $\Theta_1,...,\Theta_m\in\mathbb{C}^*$ and write
$$
\overline{n}=(n_1,...,n_m)^T,\quad N=N(\overline{n})=n_1+...+n_m.
$$
Assume that we have a sequence of simultaneous linear forms
\begin{equation}\label{LINFORMS}
L_{k,j}(\overline{n})=A_{k,0}(\overline{n})\Theta_j+A_{k,j}(\overline{n}),
\quad \overline{n}\in\mathbb{Z}_{\ge 1}^{m},
\end{equation}
$k=0,1,...m, j=1,...,m,$
where
\begin{equation}\label{ZI}
A_{k,j}=A_{k,j}(\overline{n})\in\mathbb{Z}_{\mathbb{I}},\quad k,j=0,1,...m,
\end{equation}
satisfy a determinant condition, say,
\begin{equation}\label{DET}
\Delta=
\left|\begin{matrix}
& A_{0,0}  &A_{0,1} &   ... & A_{0,m}\\
& A_{1,0}  &A_{1,1} &   ... & A_{1,m}\\
&          & ...    &       &        \\
& A_{m,0}  &A_{m,1} &   ... & A_{m,m}\\
\end{matrix} \right|
\ne 0
\end{equation}
Further, let 
$
a,b,c,d,b_i,e_i\in\mathbb{R}_{\ge 0}, a>0,                 
$
and suppose that
\begin{equation}\label{accdm}
c,c-dm>0,
\end{equation}
\begin{equation}\label{AeqN}
|A_{k,0}(\overline{n})|\le Q(\overline{n})=e^{q(N)},
\end{equation}
\begin{equation}\label{Lerin}
|L_{k,j}(\overline{n})|\le R_j(\overline{n})= e^{-r_j(\overline{n})},
\end{equation}
where
\begin{equation*}\label{qNaN}
q(N)=(aN+b\log N)g(N)+b_0N(\log N)^{1/2}+b_1N+b_2\log N+b_3,
\end{equation*}
\begin{equation*}\label{rin}
-r_j(\overline{n})=(dN-cn_j)g(N)+ e_0N(\log N)^{1/2}+e_1N+e_2\log N+e_3,
\end{equation*}
for all $k,j=0,1,...,m$.

Let the above assumptions be valid for all $N\ge N_l$, $l=1,2,3$ (where $l$ refers to case number) in our cases:\\
Case 1.  
\begin{equation*}\label{qN1}
\begin{cases}
g(N)=g_1(N):=1,\\
q(N)=q_1(N):=aN+b\log N,\\
-r_{j}(\overline{n})=-r_{j,1}(\overline{n}):=dN-cn_j+e_2\log N,
\end{cases}
\end{equation*}
and all other $b$'s and $e$'s are zero;\\
Case 2. 
\begin{equation*}\label{qN2rin2}
\begin{cases}
g(N)=g_2(N):=\log N,\quad b=0,\\
q(N)=q_2(N):=aN\log N+b_0N(\log N)^{1/2}+b_1N+b_2\log N+b_3,\\
-r_{j}(\overline{n})=-r_{j,2}(\overline{n}):=(dN-cn_j)\log N+e_0N(\log N)^{1/2}+e_1N+e_2\log N+e_3;
\end{cases}
\end{equation*}
Case 3. 
\begin{equation*}\label{qN3}
\begin{cases}
g(N)=g_3(N):=N,\\
q(N)=q_3(N):=aN^2+b_1N,\\
-r_{j}(\overline{n})=-r_{j,3}(\overline{n}):=(dN-cn_j)N+e_1N,
\end{cases}
\end{equation*}
and all other $b$'s and $e$'s are zero.

The following Theorem gives a unified result in the above three cases.

\begin{thm}\label{BAKERAXIOM}
Under the above assumptions  there exist explicit positive constants 
$F_l$ and $G_l$ not depending on $H$ such that
\begin{equation}\label{BAKERLOWER}
\left|\beta_0+\beta_1\Theta_1+...+\beta_m\Theta_m\right| > 
F_l \left(\prod_{j=1}^{m}(2mH_j)\right)^{-\frac{a}{c-dm}-\epsilon_l(H)}
\end{equation}
for all 
$\overline{\beta}=(\beta_0,\beta_1,...,\beta_m)^T \in 
\mathbb{Z}_{\mathbb{I}}^{m+1}\setminus\{\overline{0}\}$
and
\begin{equation}\label{HGREATERFk}
H=\prod_{j=1}^{m}(2mH_j)\ge G_l,\quad  H_j\ge h_j=\max\{1,|\beta_j|\}, 
\end{equation}
with an error term
$
\ \epsilon_l(H)\underset{H\to\infty}\to 0.
$
\end{thm}

In subsections \ref{subsectioncase1}-\ref{subsectioncase3} 
we consider the three cases more closely.

\subsection{Case 1}\label{subsectioncase1}

\begin{thm}\label{BAKERAXIOM1}
Denote $f=\frac{2}{c-dm}$ and
\begin{equation*}
A_1=\frac{acm}{c-dm}+B_1\log(ef),\quad B_1=\frac{ae_2m}{c-dm}+b.
\end{equation*}
Then
\begin{equation*}
F_1^{-1}=2e^{A_1},\quad
\epsilon_1(H)=B_1\frac{\log\log H}{\log H}
\end{equation*}
and
\begin{equation}\label{}
G_1=\max\{ m, N_1,e^{x_1/f}\}, \quad x_1=\max\{S_1,1\},
\end{equation}
where $S_1$ is the largest solution of the equation
\begin{equation}\label{}
S=f(e_2m\log S+dm^2+e_2m).
\end{equation}
\end{thm}

\subsection{Case 2}\label{subsectioncase2}

Before stating our results we introduce a function 
$z:\mathbb{R}\to \mathbb{R}$, 
the inverse function of the function 
$y(z) = z \log z$, $z \geq 1/e$, considered in \cite{HANCLETAL}.
\begin{lem} \label{}
\cite{HANCLETAL} The inverse function $z(y)$ of the function
$y(z)= z \log z$, $z \geq 1/e$,
is strictly increasing. 
Define $z_0(y)=y$ and $z_n(y)=\frac{y}{\log z_{n-1}}$ for $n\in\mathbb Z^+$. 
Suppose $y>e$, then 
$
z_1<z_3<\cdots <z<\cdots <z_2<z_0. 
$
Thus the inverse function may be given by the infinite nested logarithm fraction
\begin{equation*}\label{} 
z(y) =\underset{n\to\infty}{\lim} z_{n}(y)=\frac{y}{\log \frac{y}{\log \frac{y}{\log ...}}} 
\end{equation*}
for $y>e$. 
In particular,
\begin{equation}\label{loglog} 
z(y)< z_2(y) = \frac{y}{\log \frac{y}{\log y}} 
\end{equation}
for $y>e$. 
\end{lem}

\begin{thm}\label{BAKERAXIOM2}
Denote now $f=\frac{2}{c-dm}$ and
\begin{equation*}\label{}
A_2=b_0+\frac{ae_0m}{c-dm},\quad
B_2=a+b_0+b_1+\frac{ae_1m}{c-dm},
\end{equation*}
\begin{equation*}
C_2=am+b_2+\frac{a(dm^2+e_2m)}{c-dm},\quad D_2=b_0m+\frac{ae_0m^2}{c-dm},
\end{equation*}
\begin{equation*}
E_2=(a+b_0+b_1)m+b_2+b_3+\frac{a((2d+2e_0+e_1)m^2+(e_2+e_3)m)}{c-dm}.
\end{equation*}
Then
$
F_2^{-1}= 2e^{E_2}
$
and
\begin{equation}\label{THepsilon2}
\epsilon_2(H)=\xi(z,H):=
\end{equation}
$$
A_2\left(f\frac{z(f \log H)}{\log H}\right)^{1/2}+
B_2\frac{ z(f \log H)}{\log H}+
C_2\frac{\log z(f \log H)}{ \log H}+D_2\frac{(\log z(f \log H))^{1/2}}{\log H}
$$
with
\begin{equation}\label{}
G_2=\max\{ m, N_2,e^{(x_2\log x_2)/f},e^{e/f}\}, \quad x_2=\max\{S_2,1\},
\end{equation}
where $S_2$ is the largest solution of the equation
\begin{equation}\label{}
S\log S= f(e_0mS(\log S)^{1/2}+ e_1mS+(dm^2+e_2m)\log S     
\end{equation}
$$
+e_0 m^2(\log S)^{1/2}+2dm^2+2e_0m^2+e_1m^2+e_2m+e_3m).
$$
\end{thm}

In this case the estimate corresponding to (\ref{BAKERLOWER}) may be written as follows
\begin{equation}\label{BAKERGOOD}
\left|\beta_0+\beta_1\Theta_1+...+\beta_m\Theta_m\right|\ge 
\end{equation}
$$
F_2 \left(  z(f \log H)\right)^{-C_2}
H^{-\frac{a}{c-dm}-
A_2\left(f\frac{z(f \log H)}{\log H}\right)^{1/2}-
B_2\frac{ z(f \log H)}{\log H}-D_2\frac{(\log z(f \log H))^{1/2}}{\log H}}.
$$

Note, that 
\begin{equation}\label{zfloghz2}
z(f \log H)< z_2(f \log H)
\end{equation}
for $f\log H>e$ by (\ref{loglog}) and thus
\begin{equation}\label{EPSILON2xi}
\epsilon_2(H)=\xi(z,H)< \xi(z_2,H)
\end{equation}
for $f\log H>e$.
Write now
\begin{equation*}\label{}
\rho_2(x)=\frac{\log x}{\log x-\log\log x}.
\end{equation*}
Then (\ref{zfloghz2}) may further be estimated by using
\begin{equation}\label{EPSILON2xi2}
z_2(f \log H)\le \rho_2(x_0) f\left(1-\frac{\log f}{\log (f\log H)}\right)\frac{\log H}{\log\log H}
\end{equation}
valid for all
\begin{equation}\label{flogHx0}
f\log H\ge x_0\ge e^e,\quad H>e.
\end{equation}
Note, that if $0<c-dm\le 2$, then 
\begin{equation}\label{WEAKER}
z_2(f \log H)\le \rho_2(x_0) f\frac{\log H}{\log\log H}.
\end{equation}
By using the estimate (\ref{WEAKER}) we get the following corollary where
the lower bound in (\ref{USUAL}) is a generalization of what we see in the works on $E$-functions.

\begin{cor}\label{BAKERAXIOM2C}
Write $\rho=\rho_2(x_0)$.
If $0<c-dm\le 2$, $H>e$ and
\begin{equation}\label{3.20}
f\log H\ge x_0:=\max\{f\log m, f\log N_2, x_2\log x_2, e^e\},
\end{equation}
then
\begin{equation}\label{USUAL}
\left|\beta_0+\beta_1\Theta_1+...+\beta_m\Theta_m\right|\ge 
\end{equation}
$$
\frac{1}{2e^{E_2}(f\rho)^{C_2}}\left(\frac{\log\log H}{\log H}\right)^{C_2}
H^{-\frac{a}{c-dm}-
\frac{A_2f\sqrt{\rho}}{\sqrt{\log\log H}}
-\frac{B_2f\rho}{\log  \log H}
-\frac{D_2}{\log H}\sqrt{\log\left( \frac{f\rho \log H}{\log\log H}\right)} }.
$$
\end{cor}

In \cite{AKT}, $c=1, d=0$, so Corollary \ref{BAKERAXIOM2C} applies.

In most of the existing works only the terms corresponding to $A_2$ and $C_2$ are presented and
usually only a main term is given including the other terms implicitly. 
Hence in such a situation explicit dependence on the parameters, say for example on $m$, 
may become invisible.
Next we like to mention that all the methods applied to $E$-functions
seem to yield the situation where $A_2\ne 0$.
If we had $A_2=0$, then the terms with $B_2$ and $C_2$ would become more important. 
That would be the case if e.g one could find appropriate explicit Pad\'e type approximations 
instead of those produced by Siegel's lemma.

\subsection{Case 3}\label{subsectioncase3}

\begin{thm}\label{BAKERAXIOM3}
Now we have
\begin{equation*}\label{}
F_3^{-1}=2e^{B_3},\quad  
\epsilon_3(H)=A_3\frac{1}{\sqrt{\log H}},\quad G_3=\max\{ m, N_3,e\},
\end{equation*}
where general $A_3$ and $B_3$ are given in the proof section.
In the particular case, $b_1=e_1=0$, they read
\begin{equation*}
A_3=\frac{2acm}{(c-dm)^{3/2}},\quad 
B_3=\frac{acm^2(c+dm+2\sqrt{cdm})}{(c-dm)^2}.
\end{equation*}
\end{thm}

\section{Proofs}
\subsection{Proof of Theorem \ref{ZC}}
For $D\in\mathbb{Z}_{\ge 1}$, $D\not\equiv 0\pmod 4$ the ring of integers may be given by
$
\mathbb{Z}_{\mathbb{I}}=\mathbb{Z}+\mathbb{Z}(h+l\sqrt{-D})
$
with $h=0, l=1$, if $D\equiv 1$ or $2\pmod 4$ and $h=l=1/2$, if $D\equiv 3\pmod 4$.

We start with a simple principle.
First we define a lattice
\begin{equation*}\label{}
\lambda=\mathbb{Z}(1,0)+\mathbb{Z}(h,l\sqrt{D}),\quad \det\lambda=\sqrt{D}2^{-2h}
\end{equation*}
and a complex disk
\begin{equation*}\label{}
\mathcal{D}_R=\{\left. x+y(h+l\sqrt{-D})\in\mathbb{C}\right|\ x,y\in\mathbb{R},\  |x+y(h+l\sqrt{-D})|\le R\}.
\end{equation*} 
with a radius $R>0$ and a corresponding real disk
\begin{equation*}\label{}
\mathcal{C}_R=\{\left. (v,w)^T\in\mathbb{R}^2\right|\  v^2+w^2\le R^2\},
\quad \text{Vol\ }\mathcal{C}_R=\pi R^2.
\end{equation*} 
Then
\begin{equation}\label{principle}
x+y(h+l\sqrt{-D})\in\mathcal{D}_R\cap\mathbb{Z}_{\mathbb{I}}\quad\Leftrightarrow\quad
(x+yh,yl\sqrt{D})^T\in \mathcal{C}_R\cap \lambda.
\end{equation}

Next we define a lattice
\begin{equation}\label{}
\Lambda=\mathbb{Z}\overline{l}_1+...+\mathbb{Z}\overline{l}_{2m+2}\subseteq\mathbb{R}^{2m+2}
\end{equation}
generated by
\begin{equation*}\label{}
\begin{cases}\label{}
\overline{l}_1=(1,0,0,0,...,0,0)^T,\quad \overline{l}_2=(h,l\sqrt{D},0,0,...,0,0)^T,\\
\overline{l}_3=(0,0,1,0,...,0,0)^T,\quad \overline{l}_4=(0,0,h,l\sqrt{D},0,0,...,0,0)^T,\\
... \\
\overline{l}_{2m+1}=(0,0,...,0,0,1,0)^T,\quad \overline{l}_{2m+2}=(0,0,...,0,0,h,l\sqrt{D})^T.
\end{cases} 
\end{equation*} 
Immediately, $\det\Lambda=\left(\sqrt{D}2^{-2h}\right)^{m+1}$.

By using the following notations
$$
a+b(h+l\sqrt{-D})=-(z_1\Theta_1+...+z_m\Theta_m),\quad z_k=x_k+y_k(h+l\sqrt{-D}),
$$
\begin{equation*}\label{}
\quad v_k=x_k+y_kh,\quad w_k=y_kl\sqrt{D},
\quad x_k,y_k\in\mathbb{R},\quad k=0,1,...,m,
\end{equation*}
\begin{equation*}\label{}
R_0:=\left(\frac{2^{\tau}D^{1/4}}{\sqrt{\pi}}\right)^{m+1}\frac{1}{H_1\cdots H_m}
\end{equation*}
we define the following sets
\begin{equation*}\label{}
\mathcal{D}=\{\left.(z_0,z_1,...,z_m)^T\in\mathbb{C}^{m+1}\right|\
|z_0-(a+b(h+l\sqrt{-D}))|\le R_0;\ 
|z_k|\le H_k,\  k=1,...,m\},
\end{equation*} 
$$
\mathcal{C}=\{\left. (v_0,w_0,v_1,w_1,...,v_m,w_m)^T\in\mathbb{R}^{2m+2}\right|\ 
$$
\begin{equation*}\label{}
(v_0-(a+bh))^2+(w_0-bl\sqrt{D})^2\le R_0^2,\ v_k^2+w_k^2\le H_k^2,\  k=1,...,m\}.
\end{equation*} 
First we note that $\mathcal{C}$ is a symmetric convex body. 
For the volume of $\mathcal{C}$ we get
\begin{equation*}\label{}
Vol\ \mathcal{C}=\int...\int
\left(\underset{(v_0-(a+bh))^2+(w_0-bl\sqrt{D})^2\le R_0^2}{\int\int}dv_0dw_0\right)dv_1dw_1\cdots dv_mdw_m=
\end{equation*} 
\begin{equation*}\label{}
\pi R_0^2\int...\int
\left(\underset{v_1^2+w_1^2\le H_1^2}{\int\int}dv_1dw_1\right)dv_2dw_2\cdots dv_mdw_m=
...=\pi^{m+1}H_1^2\cdots H_m^2R_0^2=
\end{equation*}   
\begin{equation*}\label{}
\pi^{m+1}H_1^2\cdots H_m^2
\left(\frac{2^{2\tau}\sqrt{D}}{\pi}\right)^{m+1}\frac{1}{H_1^2\cdots H_m^2}=
2^{2m+2}\left(\frac{\sqrt{D}}{2^{2h}}\right)^{m+1}=2^{2m+2}\det\Lambda.
\end{equation*} 
Thus by Minkowski's convex body theorem, see \cite{SCHMIDT}, there exists a non-zero lattice vector
\begin{equation}\label{}
(x_0+y_0h,y_0l\sqrt{D},...,x_m+y_mh,y_ml\sqrt{D})^T\in\mathcal{C}\cap\Lambda\setminus\{\overline{0}\}.
\end{equation}
Consequently, by the above principle (\ref{principle}), we get a non-zero integer vector
\begin{equation*}\label{}
(\beta_0,\beta_1,...,\beta_m)^T=(x_0+y_0(h+l\sqrt{-D}),...,x_m+y_m(h+l\sqrt{-D}))^T
\in\mathcal{D}\cap\mathbb{Z}_{\mathbb{I}}^{m+1}\setminus\{\overline{0}\}
\end{equation*}
with
$
|\beta_k|\le H_k,\ k=1,...,m,
$
satisfying
\begin{equation}\label{linform2C}
|\beta_0+\beta_1\Theta_1+...+\beta_m\Theta_m|
\le\left(\frac{2^{\tau}D^{1/4}}{\sqrt{\pi}}\right)^{m+1}\frac{1}{H_1\cdots H_m}.\qed
\end{equation}

\subsection{Proof of Theorems \ref{BAKERAXIOM}--\ref{BAKERAXIOM3}}

Our proof starts in a classical manner and after that we give a rough description 
how to get Baker type estimates. 
Next we will introduce our tuning process which allows us to continue from the classical startup.

\subsubsection{A classical start}
We use the notation
\begin{equation*}
\Lambda:=\beta_0+\beta_1\Theta_1+...+\beta_m\Theta_m,\quad \beta_j\in\mathbb{Z}_{\mathbb{I}}
\end{equation*}
for the linear form to be estimated.
Using our simultaneous linear forms 
\begin{equation*}\label{}
L_{k,j}(\overline{n})=A_{k,0}(\overline{n})\Theta_j+A_{k,j}(\overline{n})
\end{equation*}
from (\ref{LINFORMS}) we get
\begin{equation}\label{G+R}
A_{k,0}\Lambda=\Omega_k+\beta_1 L_{k,1}(\overline{n})+...+\beta_m L_{k,m}(\overline{n}), 
\end{equation}
where
\begin{equation}\label{GKOK}
\Omega_k=\Omega_k(\overline{n})=A_{k,0}(\overline{n})\beta_0- \beta_1 A_{k,1}(\overline{n})-...
-\beta_m A_{k,m}(\overline{n})\in\mathbb{Z}_{\mathbb{I}}.
\end{equation}
If now $\Omega_k\ne 0$, then by \eqref{AeqN}, \eqref{Lerin}, \eqref{HGREATERFk},
\eqref{G+R} and \eqref{GKOK} we get
\begin{equation*}
1 \le |\Omega_k|=|A_{k,0}\Lambda-(\beta_1 L_{k,1}+...+\beta_m L_{k,m})|\le 
\end{equation*}
\begin{equation}\label{KIIKKU}
|A_{k,0}||\Lambda |+\sum_{j=1}^{m}|\beta_j||L_{k,j}|\le 
Q(\overline{n})|\Lambda |+\sum_{j=1}^{m}H_j R_j(\overline{n}). 
\end{equation}
Here we want to have, say
\begin{equation}
\sum_{j=1}^{m}H_j R_j(\overline{n})\le\frac{1}{2},
\end{equation}
in order to get a lower bound 
\begin{equation}\label{}
1\le 2|\Lambda|Q(\overline{n})
\end{equation}
for our linear form $\Lambda$.

\subsubsection{A rough version}
Here we outline a rough  version of  the proof by studying the case
$b=b_0=b_1=b_2=b_3=e_0=e_1=e_2=e_3=0$, for simplicity.
It starts by fixing the remainders and heights: 
\begin{equation}\label{ONEPAGE}
H_jR_j(\overline{n})= \frac{1}{2m}\quad\Leftrightarrow\quad 
2mH_j=e^{r_j(\overline{n})}=e^{(-dN+cn_j)g(N)}\quad\Rightarrow
\end{equation}
\begin{equation*}
e^{(-dmN+c\sum_{j=1}^{m}n_j)g(N)}=e^{(c-dm)Ng(N)}=\prod_{j=1}^{m}(2mH_j)\quad\Rightarrow
\end{equation*}
\begin{equation*}
Q(\overline{n})=e^{aNg(N)}=\left(\prod_{j=1}^{m}(2mH_j)\right)^{\frac{a}{c-dm}}\quad\Rightarrow
\end{equation*}
\begin{equation*}
1\le 2|\Lambda|Q(\overline{n})=2|\Lambda|\left(\prod_{j=1}^{m}(2mH_j)\right)^{\frac{a}{c-dm}}.
\end{equation*}

\subsubsection{Tuning}
Now a direct generalization of the second equality of \eqref{ONEPAGE} would be
\begin{equation}\label{log2mhi}
r_j(\overline{n})=\log (2mH_j),
\end{equation}
where
\begin{equation*}
r_j(\overline{n})=(-dN+cn_j)g(N)-e_0N(\log N)^{1/2}-e_1N-e_2\log N-e_3.
\end{equation*}
However,  (\ref{log2mhi}) will be too rough and thus we tune it 
into right frequency by defining
\begin{equation}\label{freq}
B_j=\log (2mH_j)+dm\hat g_l(W)+e_0m((\log W)^{1/2}+2)+e_1m+e_2,
\end{equation}
where 
\begin{equation*}
\hat g_1(W)=1,\quad \hat g_2(W)=\log W +2,\quad \hat g_3(W)=2W+m,
\end{equation*}
corresponding to our three cases.
Now we state a new system of equations
\begin{equation}\label{SsiS}
\sum_{j=1}^{m}w_j=W,
\end{equation}
\begin{equation}\label{Ssim}
r_j(\overline{w})=B_j,\quad j=1,...,m.
\end{equation}
Here (\ref{Ssim}) reads
\begin{equation}\label{pikkus}
(-dW+cw_j)g(W)-e_0W(\log W)^{1/2}-e_1W-e_2\log W-e_3=
\end{equation}
$$
\log (2mH_j)+dm\hat g_l(W)+e_0 m((\log W)^{1/2}+2)+e_1m+e_2
$$
which by (\ref{SsiS}) gives
\begin{equation}\label{TRANEQ}
(c-dm)Wg(W)-dm^2\hat g_l(W)-e_0mW(\log W)^{1/2}-e_1mW-e_2m\log W
\end{equation}
$$
-e_3m-e_0 m^2((\log W)^{1/2}+2)-e_1m^2-e_2m=\log H.
$$
The equation (\ref{TRANEQ}) has a solution $W\ge m$, if $H$ is big enough.
Then we choose the largest, say $S:=W_L\ge m$. 
(Any solution $W\ge 1$ would be satisfactory but for technical reasons we choose $W\ge m$.)
From our assumptions it follows that 
$m\ge 2$, $c>0$, $g(S)\ge 1$, $g_l(S)\ge 1$ for $l=1,2,3$, and $H_j\ge 1$ for $j=1,...,m$.
Hence $B_j\ge \log 4$ for $j=1,...,m$, which by \eqref{pikkus} implies
\begin{multline}\label{SFRAC}
s_j:=w_j=\frac{B_j+e_0S(\log S)^{1/2}+e_1S+e_2\log S+e_3+dSg(S)}{cg(S)}\\
>\frac{\log 4}{cg(S)}>0.
\end{multline}
Consequently, also the estimate \eqref{SFRAC} is valid for 
$H$ big enough (independently of each individual term $H_j$).

Put $\sigma_j=\left\lfloor s_j\right\rfloor$ and write
$\overline{\sigma}=(\sigma_1,...,\sigma_m)^T$,
$\overline{1}=(1,...,1)^T$, then
\begin{equation}\label{}
\overline{\sigma}\le \overline{s}<\overline{\sigma}+\overline{1}.
\end{equation}
First we note that
\begin{equation}\label{NSm}
T:=N(\overline{\sigma}+\overline{1})=N(\overline{\sigma})+m\le N(\overline{s})+m= S+m,\quad S<T.
\end{equation}
Next we give an estimate for the difference
\begin{multline}\label{}
r_j(\overline{s})-r_j(\overline{\sigma}+\overline{1})=\\
(-dN(\overline{s})+cs_j)g(S)-e_0N(\overline{s})(\log N(\overline{s}))^{1/2}-
e_1N(\overline{s})-e_2\log N(\overline{s})-e_3\\
-((-dN(\overline{\sigma}+\overline{1})+c(\sigma_j+1))g(N(\overline{\sigma}+\overline{1}))
-e_0N(\overline{\sigma}+\overline{1})(\log N(\overline{\sigma}+\overline{1}))^{1/2}\\
-e_1N(\overline{\sigma}+\overline{1})-e_2\log N(\overline{\sigma}+\overline{1})-e_3)=\\
d(Tg(T)-Sg(S))+c(s_jg(S)-(\sigma_j+1)g(T))\\
+e_0(T(\log T)^{1/2}-S(\log S)^{1/2})+e_1(T-S)+e_2(\log T-\log S).
\end{multline}
By $s_j<\sigma_j+1$, the increasing property of $g(x)$ and the mean value theorem we get 
\begin{multline}\label{}
r_j(\overline{s})-r_j(\overline{\sigma}+\overline{1})\le\\
d(Tg(T)-Sg(S))+e_0m((\log S)^{1/2}+2)+e_1m+e_2.
\end{multline}
Hence
\begin{equation}\label{BiMi}
r_j(\overline{s})<r_j(\overline{\sigma}+\overline{1}) +dm\hat g_l(S)+e_0m((\log S)^{1/2}+2)+e_1m+e_2,
\quad l\in\{1,2,3\},
\end{equation}
which is the reason to define (\ref{freq}).

According to the non-vanishing of the determinant (\ref{DET}) and the assumption
$\overline{\beta}=(\beta_0,\beta_1,...,\beta_m)^T \ne \overline{0}$
it follows that 
\begin{equation}
\Omega_k(\overline{\sigma}+\overline{1})\in\mathbb{Z}_{\mathbb{I}}\setminus\{0\}
\end{equation}
with some integer $k\in[0,m]$. 
Now we are ready to prove the essential estimate
\begin{equation}
\sum_{j=1}^{m} H_jR_j(\overline{\sigma}+\overline{1})= 
\sum_{j=1}^{m} H_j e^{-r_j(\overline{\sigma}+\overline{1})}\ \overset{(\ref{BiMi})}< 
\end{equation}
\begin{equation*}
\sum_{j=1}^{m} H_j e^{-B_j+dm\hat g_l(S)+e_0m((\log S)^{1/2}+2)+e_1m+e_2}=\frac{1}{2}.
\end{equation*}
Hence by (\ref{KIIKKU}) we get
\begin{equation}
1<2|\Lambda|Q(\overline{\sigma}+\overline{1})=
2|\Lambda|e^{q(N(\overline{\sigma}+\overline{1}))}\le 2|\Lambda|e^{q(S+m)},
\end{equation}
where
$$
q(S+m)=(a(S+m)+b\log(S+m))g(S+m)+
$$
\begin{equation*}
b_0(S+m)(\log (S+m))^{1/2}+b_1(S+m)+b_2\log (S+m)+b_3.
\end{equation*}
Because $g(x)$ is increasing we get
\begin{equation}
g(S+m)=g(S)+mV(S),\quad V(S)=\underset{S\le x\le S+m}\max \{g'(x)\}.
\end{equation}
Or, remembering the assumption $m\le S$, we may use the following estimates 
\begin{equation}\label{logestimates}
\log (S+m)\le \log S +1,\quad 
(\log (S+m))^{1/2}\le (\log S)^{1/2}+1.
\end{equation}
Consequently
\begin{equation}
q(S+m)\le aSg(S)+Y(S),
\end{equation}
where
$$
Y(S)=amg(S)+amSV(S)+am^2V(S)+bg(S+m)\log(S+m)+
$$
\begin{equation*}
b_0(S+m)(\log (S+m))^{1/2}+b_1(S+m)+b_2\log (S+m)+b_3.
\end{equation*}
From (\ref{TRANEQ}) we get
\begin{equation}
Sg(S)=\frac{\log H}{c-dm}+\frac{X(S)}{c-dm},
\end{equation}
where
$$
X(S)=dm^2\hat g_l(S)+e_0mS(\log S)^{1/2}+e_1mS+e_2m\log S+
$$
\begin{equation*}
e_0m^2((\log S)^{1/2}+2)+e_1m^2+e_2m+e_3m.
\end{equation*}
Hence
\begin{equation}\label{OZS428}
Q(\overline{\sigma}+\overline{1})\le H^{\frac{a}{c-dm}+Z(S)},\quad
Z(S)=\frac{1}{\log H}\left(\frac{a}{c-dm}X(S)+Y(S)\right).
\end{equation}
In the following we will consider $S$ as a variable greater than $W_L$.

\subsubsection{Case 1}
We have $\hat g_1(S)=1$ and thus
$$
Z(S)=\frac{1}{\log H}\left(\frac{a}{c-dm}(dm^2+e_2m\log S+e_2m)+
am+b\log(S+m)\right)\le
$$
\begin{equation*}
\frac{1}{\log H}\left(\frac{a(dm^2+e_2m)}{c-dm}+am+b\right)+
\frac{\log S}{\log H}\left(\frac{ae_2m}{c-dm}+b\right).
\end{equation*}
Here (\ref{TRANEQ}) reads
\begin{equation}\label{TRANEQ1}
(c-dm)W-dm^2-e_2m\log W-e_2m=\log H.
\end{equation}
Let $W_1$ denote the largest solution of the equation
\begin{equation}\label{}
(c-dm)W-dm^2-e_2m\log W-e_2m=\frac{1}{2}(c-dm)W.
\end{equation}
Hence
\begin{equation}\label{x1}
(c-dm)S-dm^2-e_2m\log S-e_2m\ge \frac{1}{2}(c-dm)W_1,\quad
\end{equation}
holds for all $S\ge x_1:=\max\{W_1,W_L,m\}$.
Further, we choose $H$ such that
\begin{equation}\label{SleflogH}
x_1\le S\le f\log H,\quad f=\frac{2}{c-dm}.
\end{equation}
Thus
\begin{equation}
Z(S)\le A_1\frac{1}{\log H}+B_1\frac{\log\log H}{\log H},
\end{equation}
where
\begin{multline*}
B_1=\frac{ae_2m}{c-dm}+b,\\
A_1=\frac{adm^2}{c-dm}+am+\frac{ae_2m}{c-dm}+b+B_1\log f=\\
\frac{acm}{c-dm}+B_1\log(ef).
\end{multline*}
Hence
\begin{equation}
1<2|\Lambda|Q(\overline{\sigma}+\overline{1})\le
|\Lambda|2e^{A_1}
H^{\frac{a}{c-dm}+B_1\frac{\log\log H}{\log H}},
\end{equation}
where
$\Lambda=\beta_0+\beta_1\Theta_1+...+\beta_m\Theta_m$
is our linear form.
This proves Theorem \ref{BAKERAXIOM1}.

\subsubsection{Case 2}
Here
\begin{equation*}
q_2(S+m)\le a(S+m)\log(S+m)+
b_0(S+m)(\log (S+m))^{1/2}+
\end{equation*}
\begin{equation*}\label{}
b_1(S+m)+
b_2\log (S+m)+
b_3
\overset{(\ref{logestimates})}\le
aS\log(S)+Y(S),
\end{equation*}
\begin{equation*}
Y(S)=b_0S(\log S)^{1/2}+(a+b_0+b_1)S+(am+b_2)\log S
\end{equation*}
\begin{equation*}
+b_0m(\log S)^{1/2}+(a+b_0+b_1)m+b_2+b_3.
\end{equation*}
From (\ref{TRANEQ}) we get
\begin{equation}
S\log S=\frac{\log H}{c-dm}+\frac{X(S)}{c-dm},
\end{equation}
where
$$
X(S)=dm^2\hat g_2(S)+e_0mS(\log S)^{1/2}+e_1mS+e_2m\log S+
$$
\begin{equation*}
e_0m^2(\log S)^{1/2}+(2e_0+e_1)m^2+e_2m+e_3m,\quad \hat g_2(S)=\log S +2.
\end{equation*}
Hence, by \eqref{OZS428},
\begin{equation*}\label{}
Z(S)=\frac{1}{\log H}
\left(A_2S(\log S)^{1/2}+B_2S+C_2\log S +D_2(\log S)^{1/2}+E\right),
\end{equation*}
where
\begin{equation*}\label{}
A_2=b_0+\frac{ae_0m}{c-dm},\quad
B_2=a+b_0+b_1+\frac{ae_1m}{c-dm},
\end{equation*}
\begin{equation*}
C_2=am+b_2+\frac{a(dm^2+e_2m)}{c-dm},\quad D=b_0m+\frac{ae_0m^2}{c-dm}.
\end{equation*}
\begin{equation*}
E_2=(a+b_0+b_1)m+b_2+b_3+\frac{a((2d+2e_0+e_1)m^2+(e_2+e_3)m)}{c-dm}.
\end{equation*}
Here (\ref{TRANEQ}) has the form
\begin{multline}\label{TRANEQ2}
(c-dm)W\log W-dm^2(\log W+2)-e_0mW(\log W)^{1/2}-e_1mW\\
-e_2m\log W-e_0 m^2((\log W)^{1/2}+2)-e_1m^2-e_2m-e_3m=\log H.
\end{multline}
Let $W_2$ denote the largest solution of the equation 
\begin{multline}\label{x2}
(c-dm)W\log W-dm^2(\log W+2)-e_0mW(\log W)^{1/2}-e_1mW\\
-e_2m\log W-e_0 m^2((\log W)^{1/2}+2)-e_1m^2-e_2m-e_3m=\frac{c-dm}{2}W\log W.
\end{multline}
Assume then  $S\ge x_2:=\max\{W_2,W_L,m\}$.
Analogously to Case 1 we may choose $H$ such that
\begin{equation}\label{SloSleflogH}
 S\log S\le f \log H,\quad f=\frac{2}{c-dm}.
\end{equation}
By (\ref{loglog}) we get 
\begin{equation}\label{SlezflogH}
S\le z(f \log H)\le z_2(f \log H)=\frac{f \log H}{\log \frac{f \log H}{\log(f\log H)}}
\end{equation}
valid for 
\begin{equation}\label{flogHe}
f\log H>e.
\end{equation}
Note the estimate
\begin{equation*}\label{}
\frac{S(\log S)^{1/2}}{\log H}=
\frac{S^{1/2}(S\log S)^{1/2}}{\log H}\le
\left(f\frac{z(f \log H)}{\log H}\right)^{1/2}\le
\left(f\frac{z_2(f \log H)}{\log H}\right)^{1/2},
\end{equation*}
too.
By using the notation
\begin{equation*}
\xi(z,H)=A_2\left(f\frac{z(f \log H)}{\log H}\right)^{1/2}+
B_2\frac{ z(f \log H)}{\log H}+
C_2\frac{\log z(f \log H)}{ \log H}+
D_2\frac{(\log z(f \log H))^{1/2}}{\log H}
\end{equation*}
we have
\begin{equation}\label{QNLE2}
Q(\overline{\sigma}+\overline{1})\le 
H^{\frac{a}{c-dm}+Z(S)}\le
e^{E_2}H^{\frac{a}{c-dm}+\xi(z,H)},
\end{equation}
where the error term satisfies
\begin{equation}\label{KSI2}
\xi(z,H)\le \xi(z_2,H).
\end{equation}

Note that
\begin{equation*}
 B_2\frac{ z(f \log H)}{\log H}=o\left(A_2\left(f\frac{z(f \log H)}{\log H}\right)^{1/2}\right)
\end{equation*}
and similarly to the terms involving $C_2$ and $D_2$. Thus 
\begin{equation*}
A_2\left(f\frac{ z(f \log H)}{\log H}\right)^{1/2}
\end{equation*}
will be the main error term, for any $H$ big enough, if $A_2\ne 0$.

Further, we note that the estimate (\ref{QNLE2}) may be written as follows
\begin{equation}\label{QNLE3}
Q(\overline{\sigma}+\overline{1})\le e^{E_2}\left(  z(f \log H)\right)^{C_2}
H^{\frac{a}{c-dm}+
A_2\left(f\frac{z(f \log H)}{\log H}\right)^{1/2}+
B_2\frac{ z(f \log H)}{\log H}+D_2\frac{(\log z(f \log H))^{1/2}}{\log H} }
\end{equation}
which by (\ref{KSI2}) implies
\begin{equation}\label{QNLE4}
Q(\overline{\sigma}+\overline{1})\le e^{E_2}\left(  z_2(f \log H)\right)^{C_2}
H^{\frac{a}{c-dm}+
A_2\left(f\frac{z_2(f \log H)}{\log H}\right)^{1/2}+
B_2\frac{ z_2(f \log H)}{\log H}+D_2\frac{(\log z_2(f \log H))^{1/2}}{\log H} }.
\end{equation}
Next we shall prove the estimates (\ref{EPSILON2xi2}), (\ref{WEAKER}) under the assumption (\ref{flogHx0}).
First we get
\begin{equation*}\label{}
z_2(y)=\frac{y}{\log y-\log\log y}
\le \frac{\log x_0}{\log x_0-\log\log x_0}\frac{y}{\log y}
=\rho_2(x_0)\frac{y}{\log y}
\end{equation*}
to be valid for all $y\ge x_0\ge e^e$.
Further, we have
\begin{equation*}\label{}
z_2(fy)\le \rho_2(x_0) f\frac{y}{\log fy}=\rho_2(x_0) f\left(1-\frac{\log f}{\log fy}\right)\frac{y}{\log y},
\end{equation*}
for $fy\ge x_0$.
In particular, we have
\begin{equation}\label{LOGE4}
z_2(f \log H)\le \rho_2(x_0) f\left(1-\frac{\log f}{\log (f\log H)}\right)\frac{\log H}{\log\log H}\le
\rho_2(x_0) f\frac{\log H}{\log\log H}
\end{equation}
for all
\begin{equation*}\label{LOGE5}
f\log H\ge x_0\ge e^e,\quad H>e,
\end{equation*}
where the last inequality in \eqref{LOGE4} is valid with $0<c-dm\le 2$.
Hence
\begin{equation}\label{QNLE5}
Q(\overline{\sigma}+\overline{1})\le e^{E_2}\left(f\rho\frac{\log H}{\log \log H}\right)^{C_2}
H^{\frac{a}{c-dm}+
\frac{A_2f\sqrt{\rho}}{\sqrt{\log\log H}}+
\frac{B_2f\rho}{\log  \log H}
+\frac{D_2}{\log H}\sqrt{\log\left( \frac{f\rho \log H}{\log\log H}\right)} },
\end{equation}
if $\rho\ge \rho_2(x_0)$, by (\ref{KSI2}), (\ref{LOGE4}).
Now substitute (\ref{QNLE2}), (\ref{QNLE3})  and (\ref{QNLE5}), respectively, into
\begin{equation}
1<2|\Lambda|Q(\overline{\sigma}+\overline{1})
\end{equation}
proving (\ref{THepsilon2}), (\ref{BAKERGOOD}) and (\ref{USUAL}).  
This ends the proof of Theorem \ref{BAKERAXIOM2} and Corollary \ref{BAKERAXIOM2C}.

\subsubsection{Case 3}
Here $\hat g_3(S)=2S+m$, so (\ref{TRANEQ}) reads
\begin{equation}\label{TRANEQ3}
(c-dm)W^2-(2dm^2+e_1m)W-dm^3-e_1m^2=\log H.
\end{equation}
Now we simply choose the larger solution  
\begin{equation}\label{}
S=\frac{2dm^2+e_1m+\sqrt{(2dm^2+e_1m)^2+4(dm^3+e_1m^2+\log H)(c-dm)}}{2(c-dm)}.
\end{equation}
For convenience, we will use the following estimate
\begin{equation}\label{SYLARAJA3}
1=S_3\le S\le v_1+v_2\sqrt{\log H}, \quad v_2=\frac{1}{\sqrt{c-dm}},
\end{equation}
\begin{equation*}
v_1= \frac{2dm^2+e_1m+\sqrt{e_1^2m^2+4cdm^3+4ce_1m^2}}{2(c-dm)}. 
\end{equation*}
Now, by using \eqref{TRANEQ3} and \eqref{SYLARAJA3}, we get
\begin{multline*}\label{}
q_3(S+m)=a(S+m)^2+b_1(S+m)=\\
\frac{a}{c-dm}\log H+\left(\frac{a(2dm^2+e_1m)}{c-dm}+2am+b_1\right)S+\frac{a(dm^3+e_1m^2)}{c-dm}+am^2+b_1m\le\\
\frac{a}{c-dm}\log H+v_2w_1\sqrt{\log H}+v_1w_1+w_2,
\end{multline*}
where
\begin{equation*}\label{}
w_1=\frac{a(2dm^2+e_1m)}{c-dm}+2am+b_1,\quad w_2=\frac{a(dm^3+e_1m^2)}{c-dm}+am^2+b_1m.
\end{equation*}
Hence
\begin{equation*}
Q(\overline{\sigma}+\overline{1})\le H^{\frac{a}{c-dm}+\frac{B_3}{\log H}+\frac{A_3}{\sqrt{\log H}}}=
e^{B_3} H^{ \frac{a}{c-dm}+A_3 \frac{1}{\sqrt{\log H}} },
\end{equation*}
where
\begin{equation*}
A_3=v_2w_1,\quad 
B_3=v_1w_1+w_2.
\end{equation*}
In particular, if $b_1=e_1=0$, then
\begin{equation*}
A_3=\frac{2acm}{(c-dm)^{3/2}},\quad 
B_3=\frac{acm^2(c+dm+2\sqrt{cdm})}{(c-dm)^2}.
\end{equation*}
This proves Theorem \ref{BAKERAXIOM3}.

\subsubsection{The term $G_l$}
Yet we need to determine terms $G_l$, $l=1,2,3$.
In each case, there are some assumptions imposed on $H$.
The determinant condition (\ref{DET}) and the conditions
$S\ge m$, (\ref{SleflogH}) and (\ref{SloSleflogH}) should be  satisfied.
So, if we put
$
f_1=x_1/f,\ f_2= (x_2\log x_2)/f,\ f_3=S_3
$
and suppose
\begin{equation}\label{Gkmax}
H\ge G_l:=\max\{ m, N_l,e^{f_l}\},
\end{equation}
then Theorem \ref{BAKERAXIOM} is proved.
Finally we note, that in Corollary \ref{BAKERAXIOM2C} we need the assumption (\ref{flogHx0}), too.
The condition (\ref{Gkmax}) applied in the Case 2 shows, in particular, that
\begin{equation}
f\log H\ge f\log G_2\ge x_2\log x_2
\end{equation}
and thus in (\ref{3.20}) we may choose 
$$
\rho=\frac{\log (x_0)}{\log (x_0)-\log\log (x_0)},
\quad x_0=\max\{f\log m, f\log N_2, x_2\log x_2, e^e\}.\qed
$$

\subsection*{Acknowledgements}
The author is deeply indebted to the anonymous referee and Anne-Maria Ernvall-Hyt\"onen for 
their useful remarks which improved the presentation of this work.
The work was supported by the Academy of Finland, grant 138522.

\end{document}